\documentclass[dvips, 10pt, a4paper]{article}
\usepackage{latexsym, amsmath}
\usepackage{amssymb}

\begin{document}
\begin{center}
{\Large{An optimal inequality on warped product semi-slant submanifolds of nearly Kaehler manifolds}}\\

\bigskip

\noindent{{Siraj Uddin, Falleh R. Al-Solamy, V.A. Khan and Cenap Ozel}}
\end{center}
\footnotetext{\it{This research is supported by the research grant UMRG 117/10AFR (University of Malaya)}}
\begin{abstract}
Non-existence of warped product semi-slant submanifolds of Kaehler manifolds was proved in [17], it is interesting to find their existence. In this paper, we prove the existence of warped product semi-slant submanifolds of nearly Kaehler manifolds by a characterization. To this end we obtain an inequality for the squared norm of second fundamental form in terms of the warping function and the slant angle. The equality case is also discussed.\\ 

\end{abstract}

\section{Introduction}
In 1969, Bishop and O'Neill [3] introduced the notion of warped product manifolds to construct examples of Riemannian manifolds with negative curvature. These manifolds are natural generalizations of Riemannian product manifolds. They defined these manifolds as: Let $(N_1, g_1)$ and $(N_2, g_2)$ be two Riemannian manifolds and $f$, a positive differentiable function on $N_1$.  Consider the product  manifold $N_1\times N_2$ with its canonical projections $\pi:N_1\times N_2 \to N_1$ and $\rho:N_1\times N_2 \to N_2$. The warped product $M=N_1\times_fN_2$ is the product manifold $N_1\times N_2$ equipped with the Riemannian structure such that
$$\|X\|^2=\|\pi_*(X)\|^2+f^2(\pi(p))\|\rho_*(X)\|^2\eqno(1.1)$$
for any tangent vector $X\in T_pM$,where $*$ is the symbol for the tangent maps. Thus we have $g=g_1+f^2g_2$. The function $f$ is called the warping function on $M$. It was proved in [3] that for any $X\in TN_1$ and $Z\in TN_2$, the following holds
$$\nabla_XZ=\nabla_ZX=(X\ln f)Z\eqno(1.2)$$
where $\nabla$ denote the Levi-Civita connection $M$. A warped product manifold $M=N_1\times _fN_2$ is said to be {\it{trivial}} if the warping function $f$ is constant. If $M=N_1\times{_{f} N_2}$ be a warped product manifold then $N_1$ is totally geodesic and $N_2$ is a totally umbilical submanifold of $M$, respectively [3]. For the survey on warped product submanifolds we refer to [8].

\parindent=8mm Nearly Kaehler manifolds are exactly the Tachibana manifolds initially studied in [18]. Nearly Kaehler manifolds form an interesting class of manifolds admitting a metric connection with parallel totally anti-symmetric torsion [1]. The best known example of a nearly Kaehler non-Kaehler manifold is $S^6$, a six dimensional sphere. It has an almost complex structure $J$ defined by the vector cross product in the space of purely imaginary Cayley numbers $\mathcal{O}$ which satisfies the nearly Kaehler structure. More general examples of nearly Kaehler manifolds are the homogeneous spaces $G/K$, where $G$ is a compact semi-simple Lie group and $K$ is the fixed point set of an automorphism of $G$ of order 3 [19]. Strict nearly Kaehler manifolds obtained a lot of consideration in 1980s due to their relation to Killing spinors. Th. Friedrich and R. Grunewald showed in [9] that a $6$-dimensional Riemannian manifold admits a Riemannian Killing spinor if
and only if it is nearly Kaehler. The only known $6$-dimensional strict nearly Kaehler manifolds are
$$S^6=G_2/SU(3).Sp(2)/SU(2)\times U(1),~~SU(3)/U(1)\times U(1),~~S^3\times S^3.$$
In fact, these are the only homogeneous nearly Kaehler manifolds in dimension six [4].

\parindent=8mm
On the other hand, slant  submanifolds of almost Hermitian manifolds were introduced by B.Y. Chen in [5] as a generalization of both holomorphic and totally real submanifolds. In [16],  N. Papaghiuc has introduced another class of submanifolds in almost Hermitian manifolds, called the semi-slant submanifolds. These submanifolds are also a natural generalization of holomorphic, totally real, slant  and CR submanifolds. 

\parindent=8mm
The idea of warped product submanifolds was given by B. Y. Chen in [7]. Since then, many mathematician extended this study for CR-warped product submanifolds of almost Hermitian as well as almost contact manifolds. Recently, B. Sahin [17] proved the non-existence of warped product semi-slant submanifolds of a Kaehler manifold. After that V.A. Khan and K.A. Khan [15] studied warped product semi-slant and warped product generic submanifolds of a nearly Kaehler manifold $\bar M$. They proved the non-existence case of warped product generic submanifolds $N\times_fN_T$ and in case of existence of warped product generic submanifolds $N_T\times_fN$, they obtained a characterization and discussed some other geometric properties, where $N_T$ and $N$ are holomorphic and a Riemannian submanifolds of $\bar M$, respectively. In case of warped product semi-slant the research problem is still open but only certain structures have the existence of warped product semi-slant submanifolds.

\parindent=8mm In this paper, first we prove the non-existence of warped product semi-slant submanifolds of the type $N_\theta\times_fN_T$ which generalizes the results obtained in [7] and [14] for the non-existence of warped product CR-submanifolds. In case of existence of warped products we obtain a geometric obstruction for second fundamental form in terms of the warping function and the slant angle which is also generalizing the results obtained for the inequality on CR-warped product submanifolds in both Kaehler and nearly Kaehler settings (see [7], [2]). 

\parindent=8mm The paper is organized as follows: In section 2, we recall some basic formulae and definitions. In section 3, we prove the non-existence of warped product semi-slant submanifolds $N_\theta\times_fN_T$ of a nearly Kaehler manifold $\bar M$, where $N_T$ and $N_\theta$ are holomorphic and proper slant submanifolds of $\bar M$. Also, the existence of warped product semi-slant submanifolds  $N_T\times_fN_\theta$ is proved by a characterization theorem. In section 4, we apply some results of section 3 and obtain an inequality for the second fundamental form. The equality case is also discussed.

\section{Preliminaries}
Let $\bar M$ be an almost Hermitian manifold with almost
complex structure $J$ and a Riemannian metric $g$ such that
$$(a)~~J^2=-I,~~~~~~~~~~(b)~~g(JX, JY) = g(X, Y)\eqno(2.1)$$
for all vector fields $X, Y$ on $\bar M$.

\parindent=8mm
Further let $T\bar M$ denote the tangent bundle of $\bar M$ and $\bar\nabla$, the covariant differential operator on $\bar M$ with respect to $g$. According to A. Gray [10], if the almost complex structure $J$ satisfies
$$(\bar\nabla_XJ)X=0\eqno(2.2)$$
\noindent for any $X\in T\bar M$, then the manifold $\bar M$ is
called a {\it{nearly Kaehler manifold}}. Equation (2.2) is equivalent to $(\bar\nabla_X J)Y+(\bar\nabla_Y J)X=0 $.

\parindent=8mm
For a submanifold $M$ of a Riemannian manifold
$\bar M$, the Gauss, Weingarten formulae are given by
$$\bar \nabla_X Y=\nabla_X Y+h(X,Y),~~~~~\bar\nabla_XN=-A_NX+\nabla^\perp_XN \eqno(2.3)$$
for all $X,Y\in TM$, where $\nabla$ is the induced
Riemannian connection on $M$, $N$ is a vector field normal to $\bar
M$, $h$ is the second fundamental form of $M$, $\nabla^\perp$ is the
normal connection in the normal bundle $T^\perp M$ and $A_N$ is the
shape operator of the second fundamental form. They are related as
$$g(A_NX,Y)=g(h(X,Y),N) \eqno(2.4)$$
where $g$ denotes the Riemannian metric on $\bar M$ as
well as the metric induced on $M$. The mean curvature vector $H$ of
$M$ is given by $H=\frac{1}{n} \sum\limits_{i=1}\limits^n h(e_i,e_i) $, where $n$ is the dimension of $M$ and
$\{e_1,e_2,\cdots,e_n\}$ is a local orthonormal frame of vector fields
on $M$. A submanifold $M$ of an almost Hermitian manifold $\bar M$ is said
to be {\it{totally umbilical}} if the second fundamental form
satisfies $h(X,Y)=g(X,Y)H,$ for all $X, Y\in TM$. The submanifold $M$ is {\it{totally geodesic}} if $h(X, Y)=0$, for all $X, Y\in TM$ and minimal if $H=0$.

\parindent=8mm A submanifold $M$ of an almost Hermitian manifold $\bar M$ is called {\it{holomorphic}} if, for any $p\in M$,  we have $ J(T_pM)=T_pM$, where $T_pM$ denotes the tangent space of $M$ at $p$. It is called {\it{totally real (or Lagrangian)}} if we have  $ J(T_pM)\subseteq T^\perp_pM$ for each $p\in M$, where $T^\perp_pM$ denotes the normal space of $M$ in $\bar M$ at $p$.

\parindent=8mm Now, let $\{e_1,\cdots, e_n\}$ be an orthonormal basis of tangent space $TM$ and $e_r$ belongs to the orthonormal basis $\{e_{n+1},\cdots e_m\}$ of the normal bundle $T^\perp M$, we put
$$h_{ij}^r=g(h(e_i, e_j), e_r)~~~~{\mbox{and}}~~~~\|h\|^2=\sum_{i,j=1}^ng(h(e_i, e_j), h(e_i, e_j)).$$

\parindent=8mm For a differentiable function $\varphi$ on $M$, the gradient $grad\varphi$ is defined by
$$g(grad\varphi, X)=X\varphi\eqno(2.5)$$
for any $X\in TM$. As a consequence, we have
$$\|grad\varphi\|^2=\sum_{i=1}^n(e_i(\varphi))^2.\eqno(2.6)$$

\parindent=8mm
For any $X\in TM$ and $N\in T^\perp M$, the transformations
$JX$ and $JN$ are decomposed into tangential, normal part as
$$JX=TX+FX,~~~~JN=BN+CN.\eqno(2.7)$$

\parindent=8mm
Now, denote by ${\cal P}_XY$ and ${\cal Q}_XY$ the tangential and normal
parts of $(\bar\nabla_XJ)Y$, i.e.,
$$(\bar\nabla_XJ)Y={\cal P}_XY+{\cal Q}_XY\eqno(2.8)$$
 for all $X, Y\in TM$. Making use of (2.7) and the Gauss and Weingarten formulae, the following
equations may easily be obtained
$${\cal P}_XY=(\bar\nabla_XT)Y-A_{FY}X-Bh(X,Y)\eqno(2.9)$$
$${\cal Q}_XY=(\bar\nabla_XF)Y+h(X, TY)-Ch(X, Y)\eqno(2.10)$$

\parindent=8mm
Similarly, for any $N\in T^\perp M$, denoting tangential and
normal parts of $(\bar\nabla_XJ)N$ by ${\cal P}_XN$ and ${\cal Q}_XN$
respectively, we obtain
$${\cal P}_XN=(\bar\nabla_XB)N+TA_NX-A_{CN}X \eqno(2.11)$$
$${\cal Q}_XN=(\bar\nabla_XC)N+h(BN, X)+FA_NX \eqno(2.12)$$
where the covariant derivative of $T,~F,~B$  and $C$ are defined
by
$$(\bar\nabla_X T)Y=\nabla_XTY-T\nabla_XY,~~~(\bar\nabla_XF)Y=\nabla^\perp_XFY-F\nabla_XY$$
$$(\bar\nabla_XB)N=\nabla_XBN-B\nabla^\perp_XN~~{\mbox{and}}~~(\bar\nabla_XC)N=\nabla^\perp_XCN-C\nabla^\perp_XN$$
for all $X,Y\in TM$ and $N\in T^\perp M$.

\parindent=8mm
It is straightforward to verify  the following properties of ${\cal P}$ and ${\cal Q}$, given by K.A. Khan and V.A. Khan in [12], which we enlist here for later use
$$(p_1)~~(i)~~~~{\cal P}_{X + Y}W={\cal P}_XW + {\cal P}_YW,
~~~~~~(ii)~~~{\cal Q}_{X + Y}W={\cal Q}_XW + {\cal Q}_YW,~~~$$
$$~(p_2)~~~(i)~~~~{\cal P}_X(Y+W)={\cal P}_XY + {\cal P}_XW,~~~~(ii)~~~~{\cal Q}_X(Y + W)={\cal Q}_XY + {\cal Q}_XW,$$
$$(p_3)~~(i)~~~~g({\cal P}_XY,~ W)=- g(Y, {\cal P}_XW),~~~(ii)~~~~g({\cal Q}_XY,~ N)=-g(Y, {\cal P}_XN),$$
$$~(p_4)~~~~{\cal P}_XJY + {\cal Q}_XJY=-J({\cal P}_XY + {\cal Q}_XY)~~~~~~~~~~~~~~~~~~~~~~~~~~~~~~~~~~~~~~~~~~~~~~~~~~~~$$
for all $X, Y, W\in TM$ and $N\in T^\perp M$.

\parindent=8mm
On a Riemannian submanifold $M$ of a nearly Kaehler manifold, by equations (2.2) and (2.8), we have
$$(a)~~{\cal P}_XY + {\cal P}_YX=0,~~(b)~~{\cal Q}_XY + {\cal Q}_YX= 0\eqno(2.13)$$
for any $X,~Y\in TM$.

\parindent=8mm For each non zero vector $X$ tangent to $M$ at $p$ the angle $\theta(X)$ between $JX$ and $T_pM$ is called the {\it{Wirtinger angle}} of $X$. Then $M$ is said to be a {\it{slant submanifold}} if the angle $\theta(X)$  is constant, which is independent of the choice of $p\in M$ and $X\in T_pM$. Holomorphic and totally real submanifolds are slant submanifolds with slant angles 0 and $\pi/2$, respectively. A slant submanifold is said to be {\it proper} if it is neither holomorphic nor totally real. More generally, a distribution $\mathcal{D}$ on $\bar M$ is called a {\it slant distribution} if the angle $\theta(X)$ between $JX$ and $\mathcal{D}_{x}$ has the same value $\theta$ for each $x\in\bar M$ and a non zero vector $X\in\mathcal{D}_{x}$. 

\parindent=8mm For a slant submanifold $M$ of an almost Hermitian manifold $\bar M$, we have
$$T^2(X)=-\cos^2\theta(X),\eqno(2.14)$$
where $\theta$ is the slant angle of $M$ in $\bar M$ (see [6]). As a consequence of the formula (2.14), we have
$$g(TX,TY)=\cos^2\theta g(X,Y),\eqno(2.15)$$
$$g(FX,FY)=\sin^2\theta g(X,Y),\eqno(2.16)$$
for any $X, Y\in TM$.
 
\parindent=8mm
 A natural generalization of CR-submanifolds in terms of slant distribution was given by N. Papaghiuc [16]. These submanifolds are known as semi-slant submanifolds. He defined these submanifolds as\\

\parindent=0mm
{\bf Definition 2.1.}  {\it{A Riemannian submanifold $M$ of an almost Hermitian manifold $\bar M$ is said to be a semi-slant submanifold if there exist two orthogonal distributions ${\mathcal{D}}_1$ and ${\mathcal{D}}_2$ such that}}
\begin{enumerate}
\item[{(i)}] $TM={\cal D}_1\oplus {\cal D}_2$
\item[{(ii)}] {\it{${\cal D}_1$ is holomorphic  i.e., $J{\cal D}_1={\cal D}_1$}}
\item[{(iii)}] {\it{${\cal D}_2$ is a slant distribution with slant angle $\theta\neq \frac{\pi}{2}.$}}
\end{enumerate}

\parindent=8mm For the integrability of the distributions involved in the above definition, V.A. Khan [13] proved the following result which we will use in our next section.\\

\parindent=0mm
{\bf Theorem 2.1.} [13] {\it{For a semi-slant submanifold $M$ of a nearly Kaehler manifold $\bar M$, the following statements are equivalent:}}
\begin{enumerate}
\item[{(i)}] {\it{The holomorphic distribution ${\cal D}_1$ on $M$ is integrable,}}
\item[{(ii)}] {\it{$h(X, TY)=h(TX, Y)$ and ${\cal Q}_XY=0$,}}
\item[{(iii)}] $h(X, TY)+h(TX, Y)=2(F\nabla_XY+Ch(X, Y))$
\end{enumerate}

\noindent{\it{for each $X, Y\in{\cal D}_1$.}}

\parindent=8mm If $\nu$ is the invariant subspace of the normal bundle $T^\perp M$ under $J$, then the normal bundle of a semi-slant submanifold can be decomposed as $T^\perp M=F{\cal D}_2\oplus\nu$, where $F{\cal D}_2$ is the normal sub bundle corresponding to the slant distribution ${\cal D}_2$ in the normal bundle $T^\perp M$.

\section{Warped product semi-slant submanifolds}
In this section, we study warped product submanifolds of a nearly Kaehler manifold $\bar M$, either in the form $N_\theta\times{_{f}N_T}$ or $N_T\times{_{f}N_\theta}$, where $N_T$ and $N_\theta$ are holomorphic and proper slant submanifolds of $\bar M$, respectively. These two types of warped products are the products in between the holomorphic and proper slant submanifolds of $\bar M$, we call such types of warped products as warped product semi-slant submanifolds in the same sense of warped product CR-submanifolds [7].\\

\noindent 
{\bf{Theorem 3.1.}} {\it Let $\bar M$ be a nearly  Kaehler manifold and $M=N_\theta\times{_{f}N_T}$ be a warped product submanifold of $\bar M$. Then  $M$ is Riemannian product of $N_T$ and $N_\theta$, where $N_T$ and  $N_\theta$ are holomorphic and proper slant submanifolds of $\bar M$, respectively.}\\

\noindent{\it Proof.} For any $X\in TN_T$ and $Z\in TN_\theta$, we have 
$$g(h(X, X), FZ)=g(\bar\nabla_XX, FZ)=g(\bar\nabla_XX, JZ)-g(\bar\nabla_XX, TZ).$$
Using the property of Riemannian metric $g$, we get
$$g(h(X, X), FZ)=-g(J\bar\nabla_XX, Z)+g(X, \bar\nabla_XTZ).$$
Then from (2.3) and the covariant derivative property of $J$, we obtain
$$g(h(X, X), FZ)=g((\bar\nabla_XJ)X, Z)-g(\bar\nabla_XJX, Z)+g(X, \nabla_XTZ).$$
Again, using the property of Riemannian metric $g$ and using (1.2), (2.2) and (2.3), we derive
$$g(h(X, X), FZ)=(TZ\ln f)\|X\|^2.\eqno(3.1)$$
Also, we have 
$$g(h(X, JX), FZ)=g(\bar\nabla_XJX, JZ)-g(\bar\nabla_XJX, TZ)$$
$$~~~~~~~~~~~~~~~~~~~~~~~=g(J\bar\nabla_XX, JZ)-g(J\bar\nabla_XX, TZ).$$
Using (2.1) and the property of Riemannian metric $g$, we get
$$g(h(X, JX), FZ)=g(\bar\nabla_XX, Z)+g(\bar\nabla_XX, T^2Z)+g(\bar\nabla_XX, FTZ).$$
Again by the property of Riemannian metric $g$ and by (1.2), (2.3) and (2.14), we obtain
$$g(h(X, JX), FZ)=-(Z\ln f)\|X\|^2+(Z\ln f)\cos^2\theta\|X\|^2+g(h(X, X), FTZ)$$
$$~~~=-(Z\ln f)\sin^2\theta\|X\|^2+g(h(X, X), FTZ).\eqno(3.2)$$
Interchanging $Z$ by $TZ$ in (3.2) and using (2.14), we derive
$$g(h(X, JX), FTZ)=-(TZ\ln f)\sin^2\theta\|X\|^2-\cos^2\theta g(h(X, X), FZ).\eqno(3.3)$$
Then from (3.1) and (3.3), we get
$$g(h(X, JX), FTZ)=-(TZ\ln f)\|X\|^2.\eqno(3.4)$$
Interchanging $X$ by $JX$ in (3.4), we obtain
$$g(h(X, JX), FTZ)=(TZ\ln f)\|X\|^2.\eqno(3.5)$$
Thus by (3.4) and (3.5), we conclude that $(TZ\ln f)\|X\|^2=0$, for any non zero vector $X\in TN_T$ and $Z\in TN_\theta$, which means that the warping function $f$ is constant on $M$ and thus the proof is complete.~$\blacksquare$

\parindent=8mm Now, we discuss the warped product semi-slant submanifolds of the type $N_T\times_fN_\theta$ of a nearly Kaehler manifold. In this case we prove the following result for later use.\\

\noindent 
{\bf{Lemma 3.1.}} {\it Let $M=N_T\times{_{f}N_\theta}$ be a warped product submanifold of a nearly  Kaehler manifold $\bar M$. Then}
\begin{enumerate}
\item[{(i)}] $g(h(X, Y), FZ)=0,$
\item[{(ii)}] $g(h(JX, Z), FZ)=(X\ln f)\|Z\|^2,$
\item[{(iii)}] $g(h(X, Z), FZ)=-(JX\ln f)\|Z\|^2,$
\end{enumerate}

\noindent {\it{for any $X, Y\in TN_T$ and $Z\in TN_\theta$.}}\\

\noindent{\it Proof.} The first part can be proved by an easy computation using (1.2), the nearly Kaehler character (2.2) and the orthogonality of two tangent spaces. For the other two parts of the lemma, we have
$$g(h(JX, Z), FZ)=g(\bar\nabla_ZJX, JZ)-g(\bar\nabla_ZJX, TZ)$$
for any $X\in TN_T$ and $Z\in TN_\theta$. Using (2.3), (1.2) and the property of Riemannian metric $g$, we obtain
$$g(h(JX, Z), FZ)=-g(JX, \bar\nabla_ZJZ)-(JX\ln f)g(Z, TZ)$$
$$~~~~~~~~~~~~~~~~~~~=-g(JX, (\bar\nabla_ZJ)Z)-g(JX, J\bar\nabla_ZZ).$$
Using (2.1), (2.2), we derive
$$g(h(JX, Z), FZ)=-g(X, \bar\nabla_ZZ)=g(\bar\nabla_ZX, Z).$$
Then by (2.3) and (1.2), we arrive at
$$g(h(JX, Z), FZ)=(X\ln f)\|Z\|^2,\eqno(3.6)$$
which is the second part of the lemma and the last part can be obtained by interchanging $X$ by $JX$ in (3.6). This completes the proof.~$\blacksquare$

\parindent=8mm If we interchange $Z$ by $TZ$ for any $Z\in TN_\theta$ in second and the third parts of the Lemma 3.1, we obtain the following relations, respectively
$$g(h(JX, TZ), FTZ)=(X\ln f)\cos^2\theta\|Z\|^2\eqno(3.7)$$
and
$$g(h(X, TZ), FTZ)=-(JX\ln f)\cos^2\theta\|Z\|^2.\eqno(3.8)$$

\parindent=8mm Now, we give the following result for later use proved by V.A. Khan and K.A. Khan [15].\\

\noindent 
{\bf{Lemma 3.2.}} {\it On a warped product semi-slant submanifold $M=N_T\times{_{f}N_\theta}$ of a nearly  Kaehler manifold $\bar M$}
$$g(h(X, TZ), FZ)=-g(h(X, Z), FTZ)=-\frac{1}{3}(X\ln f)\cos^2\theta\|Z\|^2\eqno(3.9)$$
{\it{for any $X\in TN_T$ and $Z\in TN_\theta$.}}

\parindent=8mm If we replace $X$ by $JX$ in (3.9), for any $X\in TN_T$ then we obtain
$$g(h(JX, TZ), FZ)=-g(h(JX, Z), FTZ)=-\frac{1}{3}(JX\ln f)\cos^2\theta\|Z\|^2.\eqno(3.10)$$

\parindent=8mm In the following theorem we prove the existence of warped product semi-slant submanifolds by a characterization.\\

\noindent 
{\bf{Theorem 3.2.}} \textit{Let $M$ be a semi-slant submanifold of a nearly Kaehler manifold $\bar{M}$ such that the holomorphic distribution $\mathcal{D}_1$ and slant distribution $\mathcal{D}_2$  both are integrable. Then $M$ is locally a warped product of holomorphic and proper slant submanifolds if and only if}\\
$$A_{FZ}X=-[\frac{1}{3}(X\mu)TZ+(JX\mu)Z]\eqno(3.11)$$
\textit{for any $X\in \mathcal{D}_1$ and any $Z\in \mathcal{D}_{2}$ and for a differentiable function $\mu$ on $M$ satisfying $W\mu =0$, for any $W\in\mathcal{D}_{2}$.}\\

\noindent{\it Proof.}  Let $M=N_T\times_fN_\theta$ be a warped product semi-slant submanifold of a nearly Kaehler manifold $\bar M$ such that $N_T$ and $N_\theta$ are holomorphic and proper slant submanifolds of $\bar M$, respectively. Then by Lemma 3.1 and Lemma 3.2, we get (3.11).

\parindent=8mm Conversely, if $M$ is a semi-slant submanifold with integrable holomorphic and slant distributions.  Then by Theorem 2.1 and  (2.10), we get
$$F\nabla_XY=h(X, TY)-Ch(X, Y)\eqno(3.12)$$
for any $X, Y\in \mathcal{D}_1$ . Also, by the hypothesis of the theorem, we have
$$g(h(X, Y), FZ)=g(A_{FZ}X, Y)=0,$$
for any $Z\in\mathcal{D}_{2}$. This means that $h(X, Y)\in\nu$, the invariant normal subbundle of $T^\perp M$ orthogonal to $F\mathcal{D}_{2}$.  Then the terms in the right hand side of (3.12) belongs to $\nu$ while the left hand side is in $FTM$ orthogonal to $\nu$. Hence, we conclude that $F\nabla_XY=0$, which means that $\nabla_XY\in \mathcal{D}_{1}$, for any $X, Y\in \mathcal{D}_{1}$, thus the leaves of $\mathcal{D}_{1}$  are totally geodesic  in $M$. Now, since $\mathcal{D}_{2}$ is also assumed to be integrable. We consider $N_\theta$ be a leaf of $\mathcal{D}_{2}$ and $h^\theta$ be the second fundamental form of the immersion of $N_\theta$ in $M$. Then for any $X\in\mathcal{D}_{1}$ and $Z, W\in\mathcal{D}_{2}$, we have
$$g(h^\theta(Z, W), JX)=g(\nabla_ZW, JX)=g(\bar\nabla_ZW, JX)=-g(J\bar\nabla_ZW, X).$$
Using the covariant differentiation property of $J$, we obtain
$$g(h^\theta(Z, W), JX)=g((\bar\nabla_ZJ)W, X)-g(\bar\nabla_ZJW, X).$$
Then from (2.7) and (2.8), we derive
$$g(h^\theta(Z, W), JX)=g(\mathcal{P}_ZW, X)-g(\bar\nabla_ZTW, X)-g(\bar\nabla_ZFW, X).$$
Using (2.3) and (2.4), we get
$$g(h^\theta(Z, W), JX)=g(\mathcal{P}_ZW, X)-g(h^\theta(Z, TW), X)-g(A_{FW} X, Z).$$
Then by the hypothesis of the theorem, we derive
$$g(h^\theta(Z, W), JX)=g(\mathcal{P}_ZW, X)-g(h^\theta(Z, TW), X)~~~~~~~~$$
$$~~~~~~~~~~~~~~~~~~~-\frac{1}{3}(X\mu)g(Z, TW)-(JX\mu)g(Z, W).\eqno(3.13)$$
Interchanging $X$ by $JX$ and $W$ by $TW$ and using (2.14), we arrive at
$$-g(h^\theta(Z, TW), X)=g(\mathcal{P}_ZTW, JX)+\cos^2\theta g(h^\theta(Z, W), JX)~~~~~~~~~~~~~~~$$
$$~~~~~~~~~~~~~~+\frac{1}{3}(JX\mu)\cos^2\theta g(Z,W)+(X\mu)g(Z, TW).\eqno(3.14)$$
Thus from (3.13) and (3.14), we obtain
$$\sin^2\theta g(h^\theta(Z, W), JX)=g(\mathcal{P}_ZW, X)+g(\mathcal{P}_ZTW, JX)~~~~~~~~~~~~~~~~~~~~~~~~~~~~~~~~~~~~~~~$$
$$~~~~~~~~~~~~~~~~~~~~~~~~~~~~~~~~~~~~~-(1-\frac{1}{3}\cos^2\theta)(JX\mu)g(Z,W)+\frac{2}{3}(X\mu)g(Z, TW).\eqno(3.15)$$
Now, we compute the second term of right hand side of (3.15) by using the properties of $\mathcal{P}$ as follows. By the property $(p_3)~ (i)$ and then $(p_4)$, we get
$$g(\mathcal{P}_ZTW, JX)=-g(TW, \mathcal{P}_ZJX)=g(TW, J\mathcal{P}_ZX).$$
Then using a property Riemannian metric $g$, we arrive at
$$g(\mathcal{P}_ZTW, JX)=-g(JTW, \mathcal{P}_ZX)=-g(T^2W, \mathcal{P}_ZX).$$
Using (2.14) and the property $(p_3)~(i)$, we get
$$g(\mathcal{P}_ZTW, JX)=-\cos^2\theta g(\mathcal{P}_ZW, X).\eqno(3.16)$$
Thus from (3.15) and (3.16), we derive
$$\sin^2\theta g(h^\theta(Z, W), JX)=\sin^2\theta g(\mathcal{P}_ZW, X)-(1-\frac{1}{3}\cos^2\theta)(JX\mu)g(Z,W)$$
$$+\frac{2}{3}(X\mu)g(Z, TW).\eqno(3.17)$$
By polarization identity, we obtain
$$\sin^2\theta g(h^\theta(Z, W), JX)=\sin^2\theta g(\mathcal{P}_WZ, X)-(1-\frac{1}{3}\cos^2\theta)(JX\mu)g(Z,W)$$
$$+\frac{2}{3}(X\mu)g(TZ, W).\eqno(3.18)$$
Then using the property of Riemannian metric $g$ and (2.13) (a) in (3.17) and (3.18), we get
$$ g(h^\theta(Z, W), JX)=(\frac{1}{3}\cot^2\theta-\csc^2\theta)(JX\mu)g(Z, W).$$
Finally, by (2.5) we arrive at
$$h^\theta(Z, W)=(\frac{1}{3}\cot^2\theta-\csc^2\theta)g(Z, W)grad\mu.\eqno(3.19)$$
Thus, from (3.19) we conclude that $N_\theta$ is totally umbilical in $M$ with mean curvature vector $H^\theta=(\frac{1}{3}\cot^2\theta-\csc^2\theta)grad\mu$. Now, we prove that $H^\theta$ is parallel corresponding to the normal connection $D^\theta$ of $N_\theta$ in $M$. For, this consider any $Y\in\mathcal{D}_1$ and $W\in\mathcal{D}_2$, we have
$$g(D^\theta_WH^\theta, Y)=[\frac{1}{3}\cot^2\theta-\csc^2\theta]g(D^\theta_Wgrad\mu, Y)~~~~~~~~~~~~~~~~~~~~~~~~$$
$$=[\frac{1}{3}\cot^2\theta-\csc^2\theta]g(\nabla_Wgrad\mu, Y)~~~~~~$$
$$~~~~~~~~~~~~~~~~~~=[\frac{1}{3}\cot^2\theta-\csc^2\theta]\{Wg(grad\mu, Y)-g(grad\mu, \nabla_WY)\}$$
$$~~~~~~~~~=[\frac{1}{3}\cot^2\theta-\csc^2\theta]\{W(Y\mu)-g(grad\mu, [W,Y])$$
$$-g(grad\mu, \nabla_YW)\}~~~~~~~~~~~~~~~~$$
$$~~~~~~~~~~~=[\frac{1}{3}\cot^2\theta-\csc^2\theta]\{Y(W\mu))+g(\nabla_Ygrad\mu, W)\}$$
$$=0,~~~~~~~~~~~~~~~~~~~~~~~~~~~~~~~~~~~~~~~~~~~~~~~$$
since $W\mu=0$, for all $W\in\mathcal{D}_2$ and thus $\nabla_Ygrad\mu\in\mathcal{D}_1$. This means that the mean curvature vector $H^\theta$ of $N_\theta$ is parallel. Thus the  spherical condition is satisfied, that is, $N_\theta$ is an extrinsic sphere in $M$. Hence, by a result of Hiepko [11], $M$ is a warped product manifold of $N_T$ and $N_\theta$, where $N_T$ and $N_\theta$ are integral submanifolds of $\mathcal{D}_1$ and $\mathcal{D}_2$, respectively.  This proves the theorem completely.~$\blacksquare$

\section{Inequality for warped product semi-slant submanifolds}
Now, we construct an inequality for the second fundamental form $h$ of the immersion of $M=N_T\times_fN_\theta$ into a nearly Kaehler manifold $\bar M$. We use some formulae and results from the previous sections to construct the following inequality.\\

\noindent{\bf{Theorem 4.1}} {\it Let M=$N_T\times {_{f} N_\theta}$ be a warped product semi-slant submanifold of a nearly Kaehler manifold $\bar M$ such that $N_T$ and $N_\theta$ are holomorphic and proper slant submanifolds of $\bar M$, respectively. Then} 
\begin{enumerate}
\item[{(i)}] {\it{The squared norm of the second fundamental form $h$ of $M$ satisfies}}
$$\|h\|^2 \geq 4\beta\{\csc^2\theta+\frac{1}{9}\cot^2\theta\}\|grad\ln f\|^2\eqno(4.1)$$
{\it{where $grad\ln f$ is the gradient of $\ln f$ and $2\beta$ is the dimension of $N_\theta.$}}
\item[{(ii)}] {\it{If the equality sign in (4.1) holds, then $N_T$ is totally geodesic in $\bar M$ and $N_\theta$ is a totally umbilical submanifold of $\bar M$. Moreover, $M$ is a minimal submanifold of $\bar M$.}}
\end{enumerate}

\noindent
{\it {Proof.}} Let $M=N_T\times {_{f} N_\theta}$ be a $n$-dimensional warped product semi-slant submanifold of a $m$-dimensional nearly Kaehler manifold $\bar M$ such that $N_T$ is a $2\alpha$-dimensional holomorphic submanifold and $N_\theta$ is a $2\beta$-dimensional proper slant submanifold of $\bar M$, respectively. Let us denote by ${\cal D}$ and ${\cal D}_\theta$, the tangent bundles on $N_T$ and $N_\theta$, respectively and let $\{e_1,\cdots,e_\alpha, e_{\alpha+1}=J e_1,\cdots, e_{2\alpha}=J e_\alpha\}$ and $\{ e_{2\alpha+1}=e^*_1,\cdots, e_{2\alpha+\beta} =e^*_\beta, e_{2\alpha+\beta+1}=\sec\theta Te^*_1,\cdots, e_{n}=\sec\theta Te^*_{\beta}\}$ be the local orthonormal frames of ${\cal D}$ and ${\cal D}_\theta$, respectively. Then, the orthonormal frames of $F{\cal D}_\theta$ and $\nu$ are $\{e_{n+1}=\tilde e_1=\csc\theta Fe^*_1,\cdots, \tilde e_\beta=\csc\theta Fe^*_{\beta}, \tilde e_{\beta+1}=\csc\theta\sec\theta FTe^*_1,\cdots, \tilde e_{2\beta}=\csc\theta\sec\theta FTe^*_{\beta}\}$ and $\{e_{n+2\beta+1}\cdots, e_{m}\}$, respectively, where $e_{n+2\beta+1}\cdots, e_{m}$ are orthonormal vectors in the invariant normal subbundle $\nu$ of $T^\perp M$. The dimensions of $F{\cal D}_\theta$ and $\nu$ will be $2\beta$ and $m-n-2\beta$, respectively. Then, by definition we have
$$\|h\|^2=\sum_{r=n+1}^{m}\sum_{i, j=1}^{n} g(h(e_i, e_j), e_r)^2.$$
Using assumed frames of $F{\cal D}_\theta$ and $\nu$, the above equation can be written as
$$\|h\|^2=\sum_{r=n+1}^{n+2\beta}\sum_{i, j=1}^{n} g(h(e_i, e_j), e_r)^2+\sum_{r=n+2\beta+1}^{m}\sum_{i, j=1}^{n}g(h(e_i, e_j), e_r)^2.\eqno(4.2)$$
In the first term of right hand side in (4.2) $e_r$ belongs to  $F{\cal D}_\theta$ while in the second term of right hand side $e_r$ belongs to $\nu$. We shall equate only first term of right hand side, then we get
$$\|h\|^2\geq\sum_{r=n+1}^{n+2\beta}\sum_{i, j=1}^{n} g(h(e_i, e_j), e_r)^2.$$
The above inequality for the given frame of $F{\cal D}_\theta$ will be
$$\|h\|^2\geq\sum_{r=1}^{2\beta}\sum_{i, j=1}^{n} g(h(e_i, e_j), \tilde e_r)^2.$$
Then for the orthonormal frames of ${\cal D}$ and ${\cal D}_\theta$, the above equation takes the form
$$\|h\|^2\geq\sum_{r=1}^{2\beta} \sum_{i, j=1}^{2\alpha} g(h(e_i, e_j), \tilde e_r)^2+2\sum_{r=1}^{2\beta}\sum_{i=1}^{2\alpha}\sum_{ j=1}^{2\beta} g(h(e_i, e^*_j),  \tilde e_r)^2$$
$$~~~~+\sum_{r=1}^{2\beta} \sum_{i, j=1}^{2\beta} g(h(e^*_i, e^*_j),  \tilde e_r)^2.\eqno(4.3)$$
Thus by Lemma 3.1 (i), the first term of the right hand side in (4.3) is zero and we shall compute just the next term and leave the third term, then we get
$$\|h\|^2\geq2\sum_{r=1}^{2\beta}\sum_{i=1}^{2\alpha}\sum_{j=1}^{2\beta}g(h(e_i, e^*_j),  \tilde e_r)^2.$$
Since $r, j=1,\cdots, 2\beta$, then the above inequality can be written for one summation as
$$\|h\|^2\geq2\sum_{i=1}^{2\alpha}\sum_{j=1}^{2\beta}g(h(e_i, e^*_j),  \tilde e_j)^2.$$
Using the frames of ${\cal D},~{\cal D}_\theta$ and $F{\cal D}_\theta$, we derive
$$\|h\|^2\geq2\csc^2\theta\sum_{i=1}^{\alpha}\sum_{j=1}^{\beta}g(h(e_i, e^*_j), Fe^*_j)^2~~~~~~~~~~$$
$$~~~~~~~~~+2\csc^2\theta\sec^4\theta\sum_{i=1}^{\alpha}\sum_{j=1}^{\beta}g(h(e_i, Te^*_j), FTe^*_j)^2$$
$$+2\csc^2\theta\sum_{i=1}^{\alpha}\sum_{j=1}^{\beta}g(h(J e_i, e^*_j), Fe^*_j)^2~$$
$$~~~~~~~~~~~~~~+2\csc^2\theta\sec^4\theta\sum_{i=1}^{\alpha}\sum_{j=1}^{\beta}g(h(J e_i, Te^*_j), FTe^*_j)^2~~$$
$$~~~~~~~~~+2\csc^2\theta\sec^2\theta\sum_{i=1}^{\alpha}\sum_{j=1}^{\beta}g(h(e_i, Te^*_j), Fe^*_j)^2$$
$$~~~~~~~~~~~~~~~~~+2\csc^2\theta\sec^2\theta\sum_{i=1}^{\alpha}\sum_{j=1}^{\beta}g(h(e_i, e^*_j), FTe^*_j)^2~~~~~~~$$
$$~~~~~~~~~~~~~+2\csc^2\theta\sec^2\theta\sum_{i=1}^{\alpha}\sum_{j=1}^{\beta}g(h(J e_i, Te^*_j), Fe^*_j)^2$$
$$~~~~~~~~~~~~~~~~+2\csc^2\theta\sec^2\theta\sum_{i=1}^{\alpha}\sum_{j=1}^{\beta}g(h(J e_i, e^*_j), FTe^*_j)^2.\eqno(4.4)$$
Then by Lemma 3.1 (ii)-(iii) and the relations (3.7), (3.8), (3.9) and (3.10), we obtain
$$\|h\|^2\geq4\csc^2\theta\sum_{i=1}^{\alpha}\sum_{j=1}^{\beta}\{1+\frac{1}{9}\cos^2\theta\}(e_i\ln f)^2g(e^*_j, e^*_j)^2$$
$$~~~~~~~~+4\csc^2\theta\sum_{i=1}^{\alpha}\sum_{j=1}^{\beta}\{1+\frac{1}{9}\cos^2\theta\}(J e_i\ln f)^2g(e^*_j, e^*_j)^2.$$
Thus from (2.6) the above expression will be
$$\|h\|^2\geq4\csc^2\theta\sum_{j=1}^{\beta}\{1+\frac{1}{9}\cos^2\theta\}\|grad\ln f\|^2g(e^*_j, e^*_j)^2$$
$$~~~~~~~~~~~=4\beta\{\csc^2\theta+\frac{1}{9}\cot^2\theta\}\|grad\ln f\|^2,$$
which is the inequality (4.1). If the equality holds in (4.1), then by (4.2) and (4.3), we get 
$$h({\cal D},{\cal D})=0,~~~h({\cal D}_\theta, {\cal D}_\theta)=0,~~h({\cal D}, {\cal D}_\theta)\subset F{\cal D}_\theta.\eqno(4.5)$$
Now, for any $Z, W\in{\cal D}_\theta$ and $X\in{\cal D}$, we have
$$g(h^\theta(Z, W), X)=g(\nabla_ZW, X)=-g(W, \nabla_ZX)=-(X\ln f)g(Z, W)$$
where $h^\theta$ is the second fundamental form of $N_\theta$ in $M$. Using (2.5), we obtain
$$g(h^\theta(Z, W), X)=-g(Z, W)g(grad\ln f, X),\eqno(4.6)$$
where $grad\ln f$ is the gradient of $\ln f$. Then from (4.6), we get
$$h^\theta(Z, W)=-g(Z, W)grad\ln f.\eqno(4.7)$$
Since for a warped product manifold $M=N_1\times_fN_2$, we know that $N_1$ is totally geodesic and $N_2$ is totally umbilical in $M$ [3]. Thus $N_T$ is totally geodesic submanifold in $M$ and with this fact the first condition of (4.5) implies that $N_T$ is totally geodesic in $\bar M$. Also, the second condition of (4.5) with (4.7) implies that $N_\theta$ is totally umbilical in $\bar M$. Moreover all conditions of (4.5) imply that $M$ is a minimal submanifold of $\bar M$. This completes the proof of the theorem.~$\blacksquare$

\bigskip

Siraj Uddin

\noindent Institute of Mathematical Sciences, Faculty of Science,University of Malaya, 50603 Kuala Lumpur, Malaysia

\noindent {\it E-mail}: {\tt siraj.ch@gmail.com}\\

\smallskip

Falleh Rijaullah Al-Solamy

\noindent Department of Mathematics, Faculty of Science, King Abdul Aziz University, 21589 Jeddah, Saudi Arabia

\noindent {\it E-mail}: {\tt falleh@hotmail.com}\\

\smallskip

Viqar Azam Khan

\noindent Department of Mathematics, Faculty of Science, Aligarh Muslim University, 202002 Aligarh, India

\noindent {\it E-mail}: {\tt viqarster@gmail.com}\\

\smallskip

Cenap Ozel

\noindent Department of Mathematics, Abant Izzet Baysal University, Bolu, Turkey

\noindent {\it E-mail}: {\tt cenap.ozel@gmail.com}\\

\end{document}